\documentclass[pdflatex]{sn-jnl}
\usepackage{graphicx}
\usepackage{multirow}
\usepackage{amsmath,amssymb,amsfonts}
\usepackage{amsthm}
\usepackage{mathrsfs}
\usepackage[title]{appendix}
\usepackage{xcolor}%
\usepackage{textcomp}%
\usepackage{manyfoot}%
\usepackage{booktabs}%
\usepackage{algorithm}%
\usepackage{algorithmicx}%
\usepackage{algpseudocode}%
\usepackage{listings}%
\usepackage{mathrsfs}
\usepackage{amsxtra}
\makeatletter
\renewcommand\@biblabel[1]{#1.} 
\makeatother

\usepackage{mathtools}
\usepackage{hyperref}
\hypersetup{
    colorlinks=true,		
  linkcolor=blue,		
  citecolor=orange,		
  filecolor=magenta,		
  urlcolor=cyan,		
  bookmarksdepth=4}
\newtheorem{theorem}{\bf Theorem}[section]  

\newtheorem{remark}{\bf Remark}[section]
\newtheorem{definition}{\bf Definition}[section]

\newtheorem{lem}{\bf Lemma}[section]
\begin{document}

\title[Some comparison results]{Some comparison results for one extension of the Bakry-Émery-Ricci tensor}

\author[1]{\fnm{Andrea M.} \sur{Mota}}\email{pesquisaacj@gmail.com}
\equalcont{These authors contributed equally to this work.}

\author[2]{\fnm{Cristiano S.} \sur{Silva}}\email{cristiano.silva@ufrr.br}
\equalcont{These authors contributed equally to this work.}

\author*[1]{\fnm{Juliana F.R.} \sur{Miranda}}\email{jfrmiranda@ufam.edu.br}

\affil*[1]{\orgdiv{Departamento de Matemática}, \orgname{Universidade Federal do Amazonas}, \orgaddress{\city{Manaus}, \postcode{69067-005}, \state{Amazonas}, \country{Brazil}}}

\affil[2]{\orgdiv{Departamento de Matemática}, \orgname{Universidade Federal de Roraima}, \orgaddress{\city{Boa Vista}, \postcode{69310-000}, \state{Roraima}, \country{Brazil}}}
 
\abstract{We consider generalizations of the drifted Laplacian and the Bakry-Émery-Ricci tensor, and we prove a version of the mean curvature comparison theorem. Consequently, we prove a Myers-type theorem and a generalization of the Cheeger-Gromoll splitting theorem.}



\maketitle

\begin{sloppypar}

\section{Introduction}
Comparison Geometry is based on the study of Riemannian manifolds through comparisons with the geometry of simpler spaces - the space forms. Among the seminal results are the famous Rauch's and Toponogov's theorems. The (local) Laplacian comparison theorem was without a doubt one of the first important results in comparison geometry. It states that given a complete Riemannian manifold $M$ with Ricci curvature bounded from below by a constant, the Laplacian of the distance function is bounded from above by the Laplacian of the distance function in a space form. This theorem can also be thought of as a mean curvature comparison theorem because locally the Laplacian of the distance function in $M$ can be interpreted as the mean curvature of a geodesic sphere centered on a fixed point in $M$. Its applications are often very advantageous, as revealed by two beautiful examples - Myers' theorem~\cite{myers}, and Cheeger-Gromoll splitting theorem~\cite{ChGr1}. Both theorems have significant implications for the topology and geometry of Riemannian manifolds. The former provides a bound on the diameter of a complete Riemannian manifold based on a lower bound on its Ricci curvature. The latter is a rigidity theorem whose roots are in Toponogov's theorem~\cite{top}. More precisely, a noncompact complete Riemannian manifold with nonnegative Ricci curvature splits isometrically as a Riemannian product with the real line as one factor. Not only because of its geometric character but also for its analytical nature, we can see that applications like these are quite natural in weighted Riemannian manifolds. In fact, some results in comparison geometry via the Laplacian have already been extended to the case of the drifted Laplacian operator, which is defined by means of a first-order perturbation of the Laplacian. Note that the most suitable tensor for the techniques used on weighted Riemannian manifolds is the Bakry-Émery-Ricci, introduced independently in~\cite{liche, bakem}. It is worth noting that this latter is a natural extension of the Ricci tensor. Wei and Wylie~\cite{ww} considered a lower bound on the Bakry-Émery-Ricci tensor in order to prove a comparison theorem via the drifted Laplacian. As an application, they obtained a generalized version of the Cheeger-Gromoll splitting theorem, originally due to Lichnerowicz~\cite{liche}.

By working with a symmetric, positive definite, and radially parallel (i.e. parallel in the radial direction) tensor $T$ on a weighted complete Riemannian manifold, we can obtain the natural extensions of the comparison results of~\cite{ww}. The main idea is based on some already known facts about the drifted Laplacian, and it aims at developing a suitable theory for a second-order elliptic differential operator in divergent form, so that the existing results for the Laplacian operator and the drifted Laplacian are recovered as particular cases of this operator. Here, we will call it $(\eta,T)$-divergent, namely:
\begin{equation}\label{etaT}
	\mathscr{L}f = div_{\eta}(T(\nabla f))=div(T(\nabla f))-\langle \nabla\eta,T(\nabla f)\rangle,
\end{equation}
where $f\in C^\infty(M)$.

On the one hand, associated with $\mathscr{L}$ is the tensor $R_{\eta,T} = R_{T}-\nabla(div_{\eta}T)^{\sharp}$, where $div_{\eta}T = divT-d\eta\circ T$, $R_{T}(X,Y) = tr(T\circ (Z\mapsto R(X,Z)Y))$ and $R(X,Z)Y$ denotes the Riemann curvature tensor of the metric $\langle,\rangle$. Note that the tensor $R_{\eta,T}$ is a natural extension of the Bakry-Émery-Ricci tensor. In Gomes and Miranda~\cite{gm}, some geometric properties of $R_{\eta,T}$ were established. More importantly, the following Bochner-type formula for the operator~\eqref{etaT} was proved
\begin{equation}\label{forboch}
	\frac{1}{2}\mathscr{L}(|\nabla f|^{2})=\langle\nabla(\mathscr{L}f), \nabla f\rangle + R_{\eta,T} (\nabla f, \nabla f) + \langle\nabla^{2}f, \nabla^{2}f\circ T\rangle  - \langle\nabla^{2}f,\nabla_{\nabla f}T\rangle.
\end{equation}
Observe that this formula naturally relates the analytic theory of elliptic differential operators and the geometry of the manifold.
On the other hand, Gromov~\cite{Grom} introduced the concept of generalized weighted mean curvature $H_{\eta,T}$ (for motivation see Section~\ref{sec2}). This curvature can be naturally justified when we consider the tensor $T$ to be divergence-free, i.e, $div T=0$. In this case, we prove that $\mathscr{L}r=H_{\eta,T}$, where $r$ is the distance function. When  $T$ is divergence-free the equation~\eqref{etaT} becomes
\begin{equation}\label{qdf}
	\mathscr{L}f = \square f-\langle \nabla\eta,T(\nabla f)\rangle=\square_{\eta}f.
\end{equation}
For convenience, we will refer to the operator in~\eqref{qdf} as the drifted Cheng-Yau operator. The basic tool used in the proof of the main results is the analysis of radial Ricci curvature. Indeed, if we want to get the radial curvature it is sufficient to proceed as in the definition of Riemann curvature. In Section~\ref{sec3} we will maintain the same natural approach to obtain $R_{T}$ in any radial direction.

We present below extensions of mean curvature comparison theorems obtained in~\cite{ww}. The proof of these results can be found in Section~\ref{sec4}. Throughout this paper, let $(M, \langle,\rangle)$ be an $n$-dimensional complete oriented Riemannian manifold, although completeness can be left out in Section~\ref{sec2}, $\eta\in C^{\infty}(M)$ and $T$ be a divergence-free symmetric positive definite $(1,1)$-tensor on $M$. We say that a $(1,1)$-tensor $T$ is bounded on $M$ if there exist positives $\varepsilon,\delta$ such that $\varepsilon\leq\langle T(X(t)), X(t)\rangle\leq\delta$ for all tangent vector field $X(t)$ along any minimal geodesic segment.

\begin{theorem}\label{tci}
    Let $T$ be bounded and radially parallel tensor on $M$ such that $\partial_{r}$ is an eigenvector of $T$. If 
   \begin{equation*}
		R_{\eta,T}(\partial_{r},\partial_{r})\geq \lambda, \quad \lambda\in\mathbb{R},
	\end{equation*}
	then given any minimal geodesic segment and $r_{0}\geq 0$, 
	\begin{equation}\label{TC1}
		H_{\eta,T}(r) - H_{\eta,T}(r_{0})\leq -\lambda(r-r_{0}),\quad \mbox{for}\quad r\geq r_{0}.
	\end{equation}
	Equality holds for some $r> r_{0}$ if and only if all the radial sectional curvatures are zero, $\nabla^{2}r =0$ and $\langle \nabla_{\partial_{r}}T(\nabla\eta),\partial_{r}\rangle =\lambda$, along the geodesic from $ r_{0}$ to $r$.
\end{theorem}

Henceforth, $M_{c}^{n+4k}(r)$ will denote
the simply connected space form of dimension $(n+4k)$ with constant curvature $c$, $H_{c}^{n+4k}(r)$  the mean curvature of the geodesic sphere in $M_{c}^{n+4k}(r)$ and $H_{c}$ the mean curvature of the space form of dimension $n$.

\begin{theorem}\label{tcii}Let $T$ be bounded and radially parallel tensor on $M$ such that $\partial_{r}$ is an eigenvector of $T$. Assume that
	\begin{equation*}
		R_{\eta,T}(\partial_{r},\partial_{r})\geq (n-1)c\delta.
	\end{equation*}
		
        \begin{itemize}
	\item [(a)] If  $\partial_{r}(\eta)\geq -a$, for some real constant $a \geq 0$, when $c>0$ assume $r\leq\pi/2\sqrt{c}$ then
	\begin{equation}\label{TC2}
		H_{\eta,T}(r)\leq\delta (H_{c}(r)+a)
		\end{equation}
	along that minimal geodesic segment from $p$. The equality holds if and only if the radial sectional curvatures are equal to $c$ and $\partial_{r}(\eta) =-a$. 
	\item [(b)] If $|\eta|\leq k$, when $c>0$ assume $r\leq\pi/4\sqrt{c}$ then
		\begin{equation*}
		H_{\eta,T}(r)\leq \delta H_{c}^{n+4k}(r)
		\end{equation*}
	along that minimal geodesic segment from $p$.  For the case in which $c>0$ and $r\in [\pi/4\sqrt{H},\pi/2\sqrt{H}\ ]$, we have
 
 \begin{align*}
  	H_{\eta,T}(r) &\leq \left(1+\frac{4k}{n-1}\cdot\frac{1}{\sin(2\sqrt{c}r)}\right)\delta H_{c}(r),
  \end{align*}
  along that minimal geodesic segment from $p$. In particular when $c=0$ we have
	\begin{equation}\label{tc2ii}
		H_{\eta,T}(r)\leq \delta \frac{(n+4k+1)}{r}.
	\end{equation}
	\end{itemize}
\end{theorem}

In Section~\ref{Application} as aplications we obtain a Myers’ Type Theorem and a generalization of Lichnerowicz-Cheeger-Gromoll Splitting Theorem.
 
 \begin{theorem}\label{split}
 	If $R_{\eta,T}\geq 0$ for some bounded function $\eta$ and $M$ contains a line, then $M=N^{n-1}\times \mathbb{R}$. 
 \end{theorem}

\begin{theorem}[Myers’ Type Theorem]\label{Myers} 
If $M$ has $R_{\eta,T}\geq (n-1)c\delta >0$ and $|\eta|\leq k$ then $M$ is compact and $ diam_{M} \leq \frac{\pi}{\sqrt{c}}+\frac{4k}{(n-1)\sqrt{c}}$.
\end{theorem}

\section{Generalized weighted mean curvature}\label{sec2}
Let $\Sigma$ be a hypersurface in $(M^{n},\langle,\rangle)$ and let $\nu$ be the unit normal vetor field to $\Sigma$ in $p$. Henceforth, we will denote geometric objects associated with $M$ with tilde \textquotedblleft$\sim$\textquotedblright, while those without tilde refer to $\Sigma$. We define the second fundamental form of $\Sigma$ with respect to $\nu$, as the $(0,2)$-tensor $II(X,Y)=\langle A(X),Y\rangle$, where $A(X)\coloneqq \widetilde{\nabla}_{X}\nu$ is the shape operator. Observe that $\alpha(X,Y)\coloneqq(\widetilde{\nabla}_{X}Y)^{\perp}=-II(X,Y)\nu$.

 By considering $\{e_{1},\cdots,e_{n-1},e_{n}=\nu\}$ a local orthonormal frame adapted to a neighborhood $U$ of $p$ in $\Sigma$ and $T$ a symmetric $(0,2)$-tensor on $M$, we can define a normal vector field that generalizes the mean curvature vector of $\Sigma$, $\textbf{H}\coloneqq tr\alpha$, by
\begin{equation}\label{cmg}
	\textbf{H}_{T}\coloneqq \displaystyle\sum_{i=1}^{n-1} T(e_{i},e_{i})\alpha(e_{i},e_{i}) \in{\mathfrak{X}(\Sigma)}.^{\perp}
\end{equation} 

If $\nu$ is an eigenvector of $T$ then $T(e_{i})\in T_{p}\Sigma$, so we can write \eqref{cmg} as follows
\begin{align*}
	\textbf{H}_{T}=\displaystyle\sum_{i=1}^{n-1}\alpha(T(e_{i}),e_{i}).
\end{align*}

The above definition was considered by Gomes and Miranda~\cite{gm} as well as by Roth~\cite{roth1,roth2}. By \eqref{cmg}, for all $p\in U$ we have
\begin{align*}
	\textbf{H}_{T}&= \displaystyle\sum_{i=1}^{n-1}T(e_{i},e_{i})\alpha(e_{i},e_{i})=-\displaystyle\sum_{i=1}^{n-1}T(e_{i},e_{i})
II(e_{i},e_{i})\nu  \\
	&=-\displaystyle\sum_{i=1}^{n-1}T(e_{i},e_{i})A(e_{i},e_{i})\nu =-\displaystyle\sum_{i=1}^{n-1}T(A(e_{i}),e_{i})\nu. 
\end{align*}

We will call $\textbf{H}_{T}$ the generalized mean curvature vector with respect to $\nu$, and 
\begin{equation*}
	H_{T}\coloneqq\displaystyle\sum_{i=1}^{n-1}T(A(e_{i}),e_{i})
\end{equation*}
the generalized mean curvature. Note that the sign of the function $H_{T}$ can change depending on the orientation chosen for $\nu$.

In the case that $\Sigma$ is a hypersurface in a Riemannian manifold $(M^{n},\langle,\rangle)$ equipped with a weighted volume form $dm=e^{-\eta}dM$, where $\eta\in C^{\infty}(M)$ and $dM$ is the Riemannian volume form associated with metric $\langle,\rangle$, we can consider a weighted mean curvature which naturally generalizes the mean curvature of $\Sigma$, defined by
\begin{equation}
	H_{\eta}=H-\langle\widetilde{\nabla}\eta,\nu\rangle,
\end{equation}
where $H$ is the mean curvature of $\Sigma$ with respect to $\nu$. This equation was introduced by Gromov~\cite{Grom}.

Our aim now is to define a generalized weighted mean curvature that extends the results obtained for the drifted Laplacian operator to the drifted Cheng-Yau  operator. For this purpose, first observe that if $\nu$ is a eigenvector of $T$ and $\{e_{1},\cdots,e_{n-1},e_{n}=\nu\}$ is a local orthonormal frame adapted to $\Sigma$, then the following lemma holds.

\begin{lem}\label{ecmp}
Let $\Sigma$ be a hypersurface of $M^{n}$. If $f\in C^{\infty}(M)$, then for each $p\in \Sigma$, we have
\begin{equation*}
\widetilde{\square}_{\eta}f = \square_{\eta}f + \langle  H_{T}\nu - (T(\widetilde{\nabla}\eta))^{\perp},\widetilde{\nabla}f\rangle + \widetilde{\nabla}^{2}f(\nu,T(\nu)).
\end{equation*}	

\begin{proof}  From the definition of the drifted Cheng-Yau operator, we get
	\begin{equation*}
	\widetilde{\square}_{\eta}f=\widetilde{\square}f-\langle\widetilde{\nabla} \eta,T(\widetilde{\nabla} f)\rangle=\displaystyle\sum_{i=1}^{n-1}\langle \widetilde{\nabla}_{e_{i}}\widetilde{\nabla}f, T(e_{i})\rangle + \langle \widetilde{\nabla}_{\nu}\widetilde{\nabla}f, T(\nu)\rangle - \langle\widetilde{\nabla} \eta,T(\widetilde{\nabla} f)\rangle.
	\end{equation*}
    
	Since $\widetilde{\nabla}f= \nabla f+ f_{\nu}\nu$ and $\widetilde{\nabla}\eta= \nabla \eta+ \eta_{\nu}\nu$, where $\nu(f)=f_{\nu}$ and $\nu(\eta)=\eta_{\nu}$, we have
	
	\begin{align*}
	\widetilde{\square}_{\eta}f &=\displaystyle\sum_{i=1}^{n-1}\langle \widetilde{\nabla}_{e_{i}}(\nabla f+ f_{\nu}\nu), T(e_{i})\rangle + \widetilde{\nabla}^{2}f(\nu,T(\nu))- \langle \nabla\eta+\eta_{\nu}\nu,T(\widetilde{\nabla}f)\rangle\\
	&= \displaystyle\sum_{i=1}^{n-1}\langle \widetilde{\nabla}_{e_{i}}\nabla f,T(e_{i})\rangle +  \displaystyle\sum_{i=1}^{n-1}\langle e_{i}(f_{\nu})\nu,T(e_{i})\rangle + f_{\nu}\displaystyle\sum_{i=1}^{n-1}\langle \widetilde{\nabla}_{e_{i}}\nu, T(e_{i})\rangle  \\
	&\quad+ \widetilde{\nabla}^{2}f(\nu,T(\nu))- \langle \nabla\eta,T(\nabla f)\rangle - \langle\eta_{\nu}\nu, T(\widetilde{\nabla}f)\rangle.
	\end{align*}
    
By Gauss Formula, namely $\widetilde{\nabla}_{e_{i}}\nabla f = \nabla_{e_{i}}\nabla f - II(e_{i},\nabla f)\nu$, it follows that
	\begin{align*}
		\widetilde{\square}_{\eta}f & = \displaystyle\sum_{i=1}^{n-1}\langle\nabla_{e_{i}}\nabla f,T(e_{i})\rangle  + f_{\nu}\displaystyle\sum_{i=1}^{n-1}\langle A(e_{i}),T(e_{i})\rangle+ \widetilde{\nabla}^{2}f(\nu,T(\nu)) \\
		&\quad - \langle \nabla\eta,T(\nabla f)\rangle -\langle (T(\widetilde{\nabla}\eta))^{\perp},\widetilde{\nabla}f\rangle \\
		&= \square_{\eta}f + f_{\nu}H_{T} + \widetilde{\nabla}^{2}f(\nu,T(\nu))-\langle (T(\widetilde{\nabla}\eta))^{\perp},\widetilde{\nabla}f\rangle	\\
		&=  \square_{\eta}f +  \langle H_{T}\nu,\widetilde{\nabla}f\rangle + \widetilde{\nabla}^{2}f(\nu,T(\nu))-\langle (T(\widetilde{\nabla}\eta))^{\perp},\widetilde{\nabla}f\rangle	\\
		&=\square_{\eta}f + \langle  H_{T}\nu - (T(\widetilde{\nabla}\eta))^{\perp},\widetilde{\nabla}f\rangle + \widetilde{\nabla}^{2}f(\nu,T(\nu)).
	\end{align*}
    \end{proof}
\end{lem}

 \begin{remark}
 	When $T= I$ in Lemma~\ref{ecmp} we obtain the expression of the drifted Laplacian operator on a hypersurface, and if in addition $\eta$ is a constant, the expression of the Laplacian operator on a hypersurface is also recovered.
 \end{remark}

We will denote by $\textbf{H}_{\eta,T}\coloneqq \textbf{H}_{T} - (T(\widetilde{\nabla}\eta))^{\perp}\in{\mathfrak{X}(\Sigma)}^{\perp}$ the generalized weighted mean curvature vector. This motivates the following definition:

\begin{definition}
	Let $\Sigma$ be a hypersurface in a weighted Riemannian manifold $(M,\langle,\rangle, e^{-\eta}dM)$ and let $T$ be a symmetric $(1,1)$-tensor on $M$ such that the unit normal vector field $\nu$ along $\Sigma$ is an eigenvector of $T$. We define
	\begin{equation}
		H_{\eta,T}\coloneqq H_{T}-\langle T(\widetilde{\nabla}\eta),\nu\rangle
	\end{equation}
	to be the generalized weighted mean curvature of $\Sigma$.
\end{definition}

\section{Auxiliary Lemmas}\label{sec3}
Henceforth, let $T$ be a symmetric positive definite $(1,1)$-tensor on $M$ and let $r(x)=d(x,p)$ be the distance function from a fixed point $p\in M$ which is differentiable in $M\setminus(Cut(p)\cup\{p\})$, with $|\nabla r| =1$. The level sets of $r$ are geodesic spheres on $M$ and the integral curves of radial field $\partial_{r}= \nabla r$ are geodesics, that is, $\partial_{r}|_{\gamma(t)} = \gamma'(t)$, where $\gamma$ is the unique minimal geodesic segment joining $\gamma(t)$ to $p$.

\begin{lem}\label{klem1}
	 Let $r$, $\partial_{r}$ and $T$ be defined as above. Then,
	\begin{itemize}
		\item [$(i)$]
$\nabla_{\partial_{r}}(T\circ \nabla^{2}r) + T\circ(\nabla^{2}r)^{2} = - T\circ R_{\partial_{r}}  + (\nabla_{\partial_{r}}T)\circ \nabla^{2}r.$ 
	\item [$(ii)$]
$tr((\nabla^{2}r)^{2}\circ T) + \partial_{r}(\square r) - \langle(\nabla_{\partial_{r}}T),\nabla^{2}r\rangle = - R_{T}(\partial_{r}, \partial_{r}).$ 
\end{itemize}
In particular,
\begin{equation}\label{flmtc1}
tr((\nabla^{2}r)^{2}\circ T) + \partial_{r}(\square r) = - R_{T}(\partial_{r}, \partial_{r}),
\end{equation}
whenever $T$ is a radially parallel tensor.
\begin{proof}
	Given $Y \in \mathfrak{X}(M)$, we proceed as in the definition of Riemannian curvature, that is,
\begin{equation}\label{lmt1}
	(\nabla_{Y}(T\circ \nabla^{2}r)) (\nabla r) =\nabla_{Y}T(\nabla_{\nabla r}\nabla r) - T(\nabla_{\nabla_{Y}\nabla r}\nabla r)= - T(\nabla_{\nabla_{Y}\nabla r}\nabla r).
\end{equation}

	Similarly,
	\begin{equation}\label{lmt2}
		(\nabla_{\nabla r}(T\circ \nabla^{2}r)) (Y) =\nabla_{\nabla r}T(\nabla_{Y}\nabla r) - T(\nabla_{\nabla_{\nabla r}Y}\nabla r).
	\end{equation}
    
Subtracting \eqref{lmt2} from \eqref{lmt1}, we obtain
\begin{equation}\label{lmt3}
	(\nabla_{Y}(T\circ \nabla^{2}r)) (\nabla r) - (\nabla_{\nabla r}(T\circ \nabla^{2}r)) (Y) =  - \nabla_{\nabla r}T(\nabla_{Y}\nabla r) + T(\nabla_{[\nabla r,Y]}\nabla r).
\end{equation}

Note that,
\begin{align}\label{lmt4}
	T\circ R(\nabla r,Y)\nabla r &= T(\nabla_{Y}\nabla_{\nabla r}\nabla r) - T(\nabla_{\nabla r}\nabla_{Y}\nabla r) + T(\nabla_{[\nabla r,Y]}\nabla r) \nonumber \\
	&= - T(\nabla_{\nabla r}\nabla_{Y}\nabla r) + T(\nabla_{[\nabla r,Y]}\nabla r).
\end{align}

From \eqref{lmt3} and \eqref{lmt4} we have
\begin{align}\label{lmt5}
	&(\nabla_{Y}(T\circ \nabla^{2}r)) (\nabla r) - (\nabla_{\nabla r}(T\circ \nabla^{2}r)) (Y) = \nonumber \\ &= - \nabla_{\nabla r}T(\nabla_{Y}\nabla r) + T\circ R(\nabla r,Y)\nabla r    
	+ T(\nabla_{\nabla r}\nabla_{Y}\nabla r)\nonumber \\
	&= T\circ R(\nabla r,Y)\nabla r - (\nabla_{\nabla r}T)(\nabla^{2}r(Y))\nonumber \\
	&=  T\circ R(\nabla r,Y)\nabla r - ((\nabla_{\nabla r}T)\circ \nabla^{2}r)(Y).
\end{align}

 On the other hand,
\begin{align}\label{lmtt}
	(\nabla_{Y}(T\circ \nabla^{2}r))(\nabla r)  &=\nabla_{Y}(T\circ \nabla^{2}r)(\nabla r) - (T\circ\nabla^{2}r)(\nabla^{2}r(Y)) \nonumber\\
    &= - (T\circ(\nabla^{2}r)^{2})(Y). 
\end{align}

Thus, $(i)$ follows from \eqref{lmt5} and \eqref{lmtt}. Part $(ii)$ is obtained by taking the trace from $(i)$.
\end{proof}
\end{lem}
\begin{remark}
Note that when $T=I$ in part $(i)$ of Lemma~\ref{klem1}, we have a Riccati equation, also called in this case the radial curvature equation. In part $(ii)$ we virtually obtain the Bochner formula for the Cheng-Yau operator applied to the distance function.
\end{remark}
In the next lemma we will need the following inequality for matrices
\begin{equation}\label{dmat}
	\langle I, B\rangle^{2}\leq|I_{m}|^{2}|B|^{2} = m|B|^{2}, \quad\mbox{isto é,}\quad |B|^{2}\geq\frac{1}{m}tr(B)^{2}
\end{equation}
with equality holding if and only if 
\begin{equation}\label{dmati}
	B=\frac{tr(B)}{m}I_{m},
\end{equation}
where $B=(b_{ij})_{m\times m}$ is a matrix of order $m$ and $I_{m}$ is the identity matrix.

\begin{lem}\label{ltc}
	Let $T$ be a bounded tensor on $M$ such that $\partial_{r}$ is an eigenvector of $T$. Then,
	\begin{itemize}
		\item [(i)]
		$tr((\nabla^{2}r)^{2}\circ T)\geq \displaystyle\frac{(\square r)^{2}}{(n-1)\delta}$.
		\item[(ii)]	
		$\square r=H_{T}$.
		\item[(iii)]
		$R_{T}(\partial_{r},\partial_{r}) \leq -\displaystyle\frac{H_{T}^{2}}{(n-1)\delta} - H_{T}'$,\;\mbox{whenever $T$ is a radially parallel tensor.}
	\end{itemize}
	
	Moreover, the equality holds in $(i)$ if and only if $\nabla^{2}r= \displaystyle\frac{H_{T}}{(n-1)\delta}I$, which implies that the equality also holds in $(iii)$. In all cases  $H_{T}$ denotes the mean curvature of the geodesic sphere of radius $r$ centered in $p$ with respect to the normal vector field $\nu= \partial_{r}$.
	\begin{proof}
	Let $\{e_{1},\cdots,e_{n-1}, e_{n}= \partial_{r}\}$ be a local orthonormal frame along $\gamma(t)$ which diagonalizes $T$, then
	\begin{align*}
		tr((\nabla^{2}r)^{2}\circ T) (\gamma(t)) &=  \displaystyle\sum_{i=1}^{n}\langle ((\nabla^{2}r)^{2}\circ T)(e_{i}),e_{i}\rangle\\ 
		&=  \displaystyle\sum_{i=1}^{n} \langle (\nabla^{2}r\circ T)(e_{i}),\nabla^{2}r(e_{i})\rangle \\
		&= \displaystyle\sum_{i,j=1}^{n}\langle  (\nabla^{2}r\circ T)(e_{i}), e_{j}\rangle\;\langle\nabla^{2}r(e_{i}),e_{j}\rangle \\
		&= \displaystyle\sum_{i,j=1}^{n}\frac{1}{\lambda_{i}}\langle  (\nabla^{2}r\circ T)(e_{i}), e_{j}\rangle\;\langle\nabla^{2}r(\lambda_{i}e_{i}),e_{j}\rangle.
	\end{align*}
    Since $T(e_{i}) = \lambda_{i}e_{i}$, for all $i = 1,\cdots, n$, we have
	\begin{equation*}
		tr((\nabla^{2}r)^{2}\circ T)(\gamma(t)) =  \displaystyle\sum_{i,j=1}^{n}\frac{1}{\lambda_{i}}\langle  (\nabla^{2}r\circ T)(e_{i}), e_{j}\rangle\;\langle(\nabla^{2}r\circ T)(e_{i}),e_{j}\rangle.
	\end{equation*}
	Using the fact that $\partial_{r}$ is an eigenvector of $T$ and $\nabla_{\partial_{r}}\partial_{r} = 0$, we get
	\begin{align}\label{ltc1}
		tr((\nabla^{2}r)^{2}\circ T)(\gamma(t)) &= \displaystyle\sum_{i,j=1}^{n-1}\left(\frac{\langle(\nabla^{2}r\circ T)(e_{i}),e_{j}\rangle}{\sqrt{\lambda_{i}}}\right)^{2} \nonumber\\
		&\stackrel{\eqref{dmat}}\geq \frac{1}{n-1}\left(\displaystyle\sum_{i=1}^{n-1}\frac{\langle(\nabla^{2}r\circ T)(e_{i}),e_{i}\rangle}{\sqrt{\lambda_{i}}}\right)^{2}.
	\end{align}
	Since $\lambda_{i} = \langle T(e_{i}),e_{i}\rangle\leq \delta$, for each $i=1,\cdots,n$, it follows that
	\begin{align}
		tr((\nabla^{2}r)^{2}\circ T) (\gamma(t))&\geq	\frac{1}{n-1}\left(\displaystyle\sum_{i=1}^{n-1}\frac{\langle(\nabla^{2}r\circ T)(e_{i}),e_{i}\rangle}{\sqrt{\delta}}\right)^{2} \nonumber\\
		&= \frac{1}{(n-1)\delta}\left(\displaystyle\sum_{i=1}^{n}\langle(\nabla^{2}r\circ T)(e_{i}),e_{i}\rangle\right)^{2} \nonumber\\
		&=\frac{(\square r)^{2}}{(n-1)\delta},\nonumber
	\end{align}
	which proves $(i)$. For $(ii)$, we have
	\begin{align*}
		\square r &= \displaystyle\sum_{i=1}^{n}\langle\nabla_{e_{i}}\nabla r, T(e_{i})\rangle = \displaystyle\sum_{i=1}^{n-1}\langle\nabla_{e_{i}}\nabla r, T(e_{i})\rangle  = \displaystyle\sum_{i=1}^{n-1}\langle A(e_{i}), T(e_{i})\rangle = H_{T}.
	\end{align*}
    
	Part $(iii)$ is an immediate consequence of $(i)$, $(ii)$ and \eqref{flmtc1}. Now, the equality holds in $(i)$ if and only if
	\begin{equation*}
		\lambda_{i} = \delta\quad\mbox{and}\quad\displaystyle\frac{\langle(\nabla^{2}r\circ T)(e_{i}),e_{i}\rangle}{\sqrt{\delta}} \stackrel{\eqref{dmati}} = \psi,
	\end{equation*}
	for each $i =1,\cdots,n-1$.
	
Note that $(\nabla^{2}r\circ T)/\sqrt{\delta}$ is the operator whose entries of its associated matrix in the basis $\{e_{1},\cdots,e_{n-1}, e_{n}= \partial_{r}\}$ are exactly $\langle(\nabla^{2}r\circ T)(e_{i}),e_{i}\rangle/\sqrt{\delta}$  and that
	\begin{equation*}
		\psi=\displaystyle\frac{\square r}{(n-1)\sqrt{\delta}}.
	\end{equation*}
	
	Since  $T(e_{i})=\delta e_{i}$ for all $i=1,\cdots,n-1$, the equality occurs in \eqref{ltc1}  if and only if 
	\begin{align*}
		\nabla^{2}r(e_{i},e_{i}) &= \psi\frac{\sqrt{\delta}}{\delta}
		=\displaystyle\frac{\square r}{(n-1)\delta}
		=\displaystyle\frac{H_{T}}{(n-1)\delta}.
	\end{align*}
\end{proof}
\end{lem}

\section{Proof of Theorems \ref{tci} and \ref{tcii}}\label{sec4}%

To prove our version of the mean curvature comparison theorem, we follow the techniques of Wei and Wylie~\cite{ww}.

\begin{proof}[\textbf{Proof of Theorem~\ref{tci}}]
    From Lemma~\ref{ltc}
\begin{equation}\label{tci2}
	\frac{H_{T}^{2}}{(n-1)\delta} + H_{T}' + R_{T}(\partial_{r},\partial_{r})\leq 0,
\end{equation}
and since $H_{\eta,T}'=H_{T}'-\langle \nabla_{\partial_{r}}T(\nabla\eta),\partial_{r}\rangle$, we obtain 
\begin{equation}\label{tci4}
	H_{\eta,T}'\leq - \frac{H_{T}^{2}}{(n-1)\delta} - R_{\eta,T}(\partial_{r},\partial_{r}).
\end{equation}
    By hypothesis on $R_{\eta,T}$, we get
	\begin{equation}\label{tci5}
		H_{\eta,T}'\leq - \frac{H_{T}^{2}}{(n-1)\delta} - \lambda,
	\end{equation}
	whence

\begin{equation}\label{tci6}
	H_{\eta,T}'\leq - \lambda.
\end{equation}

    Integrating \eqref{tci6} from $r_{0}$ to $r$ with $r_{0}\leq r$, i.e.,

\begin{equation*}
	\int_{r_{0}}^{r} (H'_{\eta,T}(t)+\lambda)\;dt\leq0,
\end{equation*}
	 we obtain the required inequality \eqref{TC1}.

Equality in \eqref{TC1} holds if and only if $H_{\eta,T}'=- \lambda$ on an interval $[r_{0},r]$. Consequently, we have from \eqref{tci5} that $H_{T}=0$ and $\langle \nabla_{\partial_{r}}T(\nabla\eta),\partial_{r}\rangle =\lambda$. By the hypothesis on $R_{\eta,T}$ and \eqref{tci4}, we obtain $R_{\eta,T} = \lambda$ and therefore $R_{T}(\partial_{r},\partial_{r})=0$.   Moreover, by the fact that $H_{T} =0$ and $R_{T}(\partial_{r},\partial_{r})=0$, it follows from \eqref{flmtc1} that $tr((\nabla^{2}r)^{2}\circ T)=0$ and from part $(i)$ of Lemma~\ref{ltc} that  $\nabla^{2}r =0$.

Now, choosing a local orthonormal frame along $\gamma(t)$ which diagonalizes $T$, from part $(i)$ of Lemma \ref{klem1} we have
    \begin{equation*}
		T\circ R(\partial_{r},e_{i})\partial_{r} =   0,
	\end{equation*}
    and since $T(e_{i})=\delta e_{i}$, for all $i=1,\cdots,n-1$, we get
	\begin{equation*}
		0 = \langle T\circ R(\partial_{r},e_{i})\partial_{r}, e_{i}\rangle = \langle R(\partial_{r},e_{i})\partial_{r}, T(e_{i})\rangle = \delta K(\partial_{r},e_{i}).
	\end{equation*}
	Finally, because $\delta$ is positive, we conclude that $K(\partial_{r},e_{i})=0$.
\end{proof}

\begin{proof}[\textbf{Proof of Theorem~\ref{tcii}}]
$(a)$ Multiplying \eqref{tci2} by $\displaystyle\frac{1}{\delta}$ and using the hypothesis on $R_{\eta,T}$, we have
\begin{align}\label{tcii1}
	\frac{H_{T}'}{\delta}&\leq - \frac{H_{T}^{2}}{(n-1)\delta^{2}} - \frac{1}{\delta}R_{T}(\partial_{r}, \partial_{r}) \nonumber \\
	&= - \frac{H_{T}^{2}}{(n-1)\delta^{2}} - \frac{1}{\delta}R_{\eta,T}(\partial_{r}, \partial_{r}) +\frac{1}{\delta}\langle\nabla_{\partial_{r}}T(\nabla\eta),\partial_{r}\rangle \nonumber \\
	&\leq - \frac{H_{T}^{2}}{(n-1)\delta^{2}} - (n-1)c +\frac{1}{\delta}\langle\nabla_{\partial_{r}}T(\nabla\eta),\partial_{r}\rangle.	
\end{align}

Applying the classical Bochner formula~\eqref{flmtc1} to the distance function in $M^{n}_{c}$, we obtain
\begin{equation}\label{tcii2}
	H_{c}'= - \frac{H_{c}^{2}}{n-1} - (n-1)c.
\end{equation}

Subtracting \eqref{tcii2} from \eqref{tcii1} we get
\begin{equation}\label{tcii3}
	\left(\frac{H_{T}}{\delta} - H_{c}\right)'\leq - \left(\frac{H_{T}^{2}}{(n-1)\delta^{2}} - \frac{H_{c}^{2}}{n-1}\right) + \frac{1}{\delta}\langle\nabla_{\partial_{r}}T(\nabla\eta),\partial_{r}\rangle.
\end{equation}

To simplify the above inequality \eqref{tcii3}, consider the function
\[
sn_{c}(r) = 
\begin{cases}
	\displaystyle\frac{1}{\sqrt{c}}sin(\sqrt{c}r), &\text{if $c>0$,} \\
	r, &\text{if $c=0,$}\\
	\displaystyle\frac{1}{\sqrt{c}}sinh(\sqrt{|c|}r), &\text{if $c<0$}
\end{cases}
\]	
which is the solution to the following initial value problem
\begin{equation*}
	sn_{c}''+ csn_{c}=0 \quad\mbox{such that}\quad sn_{c}(0)=0\quad\mbox{and}\quad sn_{c}'(0)=1.
\end{equation*}

From \eqref{tcii3} and since
\begin{equation}\label{tcii4}
	H_{c}=(n-1)\frac{sn_{c}'}{sn_{c}},
\end{equation}
 we compute
\begin{align}\label{tcii5}
	\left[sn_{c}^{2}\left(\frac{H_{T}}{\delta} - H_{c}\right)\right]'&=2sn_{c}sn_{c}'\left(\frac{H_{T}}{\delta} - H_{c}\right) +  sn_{c}^{2}\left(\frac{H_{T}}{\delta} - H_{c}\right)' \nonumber \\
	&\leq \frac{2sn_{c}^{2}H_{c}}{n-1}\left(\frac{H_{T}}{\delta} - H_{c}\right) - sn_{c}^{2}\left(\frac{H_{T}^{2}}{(n-1)\delta^{2}} - \frac{H_{c}^{2}}{n-1}\right) \nonumber \\
	&\quad +\frac{sn_{c}^{2}}{\delta}\langle\nabla_{\partial_{r}}T(\nabla\eta),\partial_{r}\rangle  \nonumber \\
	&= sn_{c}^{2}\left[-\frac{1}{n-1}\left(\frac{H_{T}}{\delta}- H_{c}\right)^{2} +  \frac{1}{\delta}\langle\nabla_{\partial_{r}}T(\nabla\eta),\partial_{r}\rangle\right] \nonumber \\
	&\quad \leq \frac{sn_{c}^{2}}{\delta}\langle\nabla_{\partial_{r}}T(\nabla\eta),\partial_{r}\rangle.
\end{align}

Integrating \eqref{tcii5} from $0$ to $r$ yields
\begin{equation}\label{tcii6}\frac{1}{\delta}sn_{c}^{2}(r)H_{T}(r) - sn_{c}^{2}(r) H_{c}(r)\leq \frac{1}{\delta}\int_{0}^{r}sn_{c}^{2}(t)\langle\nabla_{\partial_{t}}T(\nabla\eta),\partial_{t}\rangle\;dt.
\end{equation}

As $\langle\nabla_{\partial_{t}}T(\nabla\eta),\partial_{t}\rangle = \partial_{t} \langle T(\nabla\eta),\partial_{t}\rangle$, the integration by parts in~\eqref{tcii6} results in
\begin{equation}\label{tcii7}
	\int_{0}^{r}sn_{c}^{2}(t)\langle\nabla_{\partial_{t}}T(\nabla\eta),\partial_{t}\rangle\;dt = sn_{c}^{2}(r)\langle T(\nabla\eta),\partial_{r}\rangle - \int_{0}^{r}(sn_{c}^{2}(t))'\langle T(\nabla\eta),\partial_{t}\rangle\;dt.
\end{equation} 

Combining \eqref{tcii6} with \eqref{tcii7}, we have
\begin{align*}
	\frac{1}{\delta}sn_{c}^{2}(r)H_{\eta,T}(r) - sn_{c}^{2}(r) H_{c}(r) &\leq - \frac{1}{\delta}\int_{0}^{r}(sn_{c}^{2}(t))'\langle T(\nabla\eta),\partial_{t}\rangle\;dt \\
	&\leq - \frac{1}{\delta}\int_{0}^{r}(sn_{c}^{2}(t))'\langle \nabla\eta,T(\partial_{t})\rangle\;dt.
\end{align*}

Making use of the hypothesis that $T(\partial_{t}) = \mu\partial_{t}$ and $\mu = \langle T(\partial_{t}),\partial_{t}\rangle\leq\delta$ for some function $\mu\in C^{0}(M)$ we can rewrite the previous inequality as 
\begin{align}\label{tcii8}
	\frac{1}{\delta}sn_{c}^{2}(r)H_{\eta,T}(r) - sn_{c}^{2}(r) H_{c}(r) &\leq -\frac{1}{\delta}\int_{0}^{r}(sn_{c}^{2}(t))'\mu\;\partial_{t}(\eta)(t)\;dt \nonumber\\
	&\leq -\int_{0}^{r}(sn_{c}^{2}(t))'\;\partial_{t}(\eta)(t)\;dt.		
\end{align}

Since $(sn_{c}^{2}(t))' = 2sn_{c}'(t)sn_{c}(t)\geq 0$ and $\partial_{t}(\eta)(t)\geq -a$, it follows that
\begin{equation*}
	\frac{1}{\delta}sn_{c}^{2}(r)H_{\eta,T}(r) - sn_{c}^{2}(r) H_{c}(r) \leq  a\int_{0}^{r}(sn_{c}^{2}(t))'\;dt = asn_{c}^{2}(r),
\end{equation*}
and hence
\begin{equation*}
	H_{\eta,T}(r)\leq\delta (H_{c}(r)+ a).
\end{equation*}

To discuss when the equality holds, suppose that $H_{\eta,T}(r) =\delta	H_{c}(r)+ \delta a$ for some $r$ and $\partial_{t}(\eta)\geq -a$. Then from \eqref{tcii8} we see
\begin{equation*}
	\frac{1}{\delta}sn_{c}^{2}(r)(H_{\eta,T}(r) - \delta H_{c}(r))\leq -\int_{0}^{r}(sn_{c}^{2}(t))'\;\partial_{t}(\eta)(t)\;dt,
\end{equation*}
that is, 
\begin{equation}\label{tcii9}
	asn_{c}^{2}(r)\leq -\int_{0}^{r}(sn_{c}^{2}(t))'\;\partial_{t}(\eta)(t)\;dt.
\end{equation}

On the other hand,
\begin{equation}\label{tcii10}
	-\int_{0}^{r}(sn_{c}^{2}(t))'\;\partial_{t}(\eta)(t)\;dt\leq asn_{c}^{2}(r).
\end{equation}

From \eqref{tcii9} and \eqref{tcii10}, we obtain
\begin{equation*}
	asn_{c}^{2}(r)\leq -\int_{0}^{r}(sn_{c}^{2}(t))'\;\partial_{t}(\eta)(t)\;dt \leq asn_{c}^{2}(r),
\end{equation*}
and so $\partial_{t}(\eta) = -a$. By definition $H_{\eta, T} = H_{T} - \langle T(\nabla\eta),\partial_{r}\rangle$, which implies that
\begin{equation*}
	H_{\eta, T}(r)= H_{T}(r) - \mu\partial_{r}(\eta).	
\end{equation*}

Notice that equality \eqref{TC2} occurs when the equality in part $(iii)$ of Lemma~\ref{ltc} holds together with $\mu=\delta$. Thus, $H_{T}(r) = H_{\eta, T}(r) - \delta a = \delta H_{c}(r)$ and
\begin{equation*}
	\nabla^{2}r = \frac{H_{T}(r)}{(n-1)\delta}I = \frac{H_{c}(r)}{(n-1)}I \stackrel{\eqref{tcii4}}= \frac{sn_{c}'(r)}{sn_{c}(r)} I.
\end{equation*}

Thus, for all unit tangent vector $u\in T_{\gamma(t)}M$ with $u\perp\partial_{r}$, the radial sectional curvature is computed as follows
\begin{align*}
	-R(\partial_{r},u)\partial_{r} &= \left(\nabla_{\partial_{r}}\left(\frac{sn_{c}'(r)}{sn_{c}(r)}I\right)\right) u + \left(\frac{sn_{c}'(r)}{sn_{c}(r)}\right)^{2}u \\
	&=\left(\frac{sn_{c}'(r)}{sn_{c}(r)}\right)'u  +\left(\frac{sn_{c}'(r)}{sn_{c}(r)}\right)^{2}u  \\
	&=\left[\left(\frac{sn_{c}'(r)}{sn_{c}(r)}\right)'+ \left(\frac{sn_{c}'(r)}{sn_{c}(r)}\right)^{2} \right]u \\
	&=\displaystyle\frac{sn_{c}''(r)}{sn_{c}(r)} u,
\end{align*}
and we obtain
\begin{equation*}
	K(\partial_{r},u) = \langle R(\partial_{r},u)\partial_{r}, u\rangle = - \frac{sn_{c}''(r)}{sn_{c}(r)}=c.
\end{equation*}
$(b)$ Using integration by parts we have 
	\begin{equation*}
		\int_{0}^{r}(sn_{c}^{2}(t))'\;\partial_{t}(\eta)(t)\;dt= (sn_{c}^{2}(r))'\eta(r) -\int_{0}^{r}\eta(t)(sn_{c}^{2})''(t)\;dt,
	\end{equation*}
whence \eqref{tcii8} can be written as

\begin{equation}\label{tcii15}
	\frac{1}{\delta}sn_{c}^{2}(r)H_{\eta,T}(r) \leq  sn_{c}^{2}(r) H_{c}(r) - \eta(r)(sn_{c}^{2}(r))' +	\int_{0}^{r}\eta(t)(sn_{c}^{2})''(t)\;dt.
\end{equation}
Now if $|\eta|\leq k$ and $r\in (0,\pi/4\sqrt{c}\ ]$ when $c>0$, then $(sn_{c}^{2})''(t)>0$ and
\begin{equation*}
		\frac{1}{\delta}sn_{c}^{2}(r)H_{\eta,T}(r) \leq sn_{c}^{2}(r) H_{c}(r) + 2k(sn_{c}^{2}(r))'.
\end{equation*}
From \eqref{tcii4} we can see that
\begin{equation*}
	(sn_{c}^{2}(r))'= 2sn_{c}'(r)sn_{c}(r)= \frac{2}{n-1}H_{c}(r)sn_{c}^{2}(r),
\end{equation*}
so we have
\begin{equation*}
	H_{\eta,T}(r)\leq \left(1+\frac{4k}{n-1}\right)\delta H_{c}(r)= \delta H_{c}^{n+4k}(r).
\end{equation*}
In particular, when $c=0$ we have $sn_{c}(r)=r$ and therefore
\begin{equation*}
		H_{\eta,T}(r)\leq \left(1+\frac{4k}{n-1}\right)\delta \frac{(n-1)}{r}= \delta \frac{(n+4k-1)}{r}.
\end{equation*}

Now, we give an estimate for $H_{\eta,T}$ that will be used later to prove one of our applications. As in~\cite{ww}, when $c>0$ and $r\in [\pi/4\sqrt{H},\pi/2\sqrt{H}\ ]$,
 \begin{equation*}
  	\int_{0}^{r}f(t)(sn_{c}^{2})''(t)\;dt\leq k\left(\frac{2}{\sqrt{c}}-sn_{c}(2r)\right).
  \end{equation*}.

Replacing the estimate above in~\eqref{tcii15}, we have
   \begin{align}\label{est1}
 	H_{\eta,T}(r) &\leq \delta H_{c}(r) - \frac{\delta\eta(r)(sn_{c}^{2}(r))'}{sn_{c}^{2}(r)} +	\frac{\delta k}{sn_{c}^{2}(r)}\left(\frac{2}{\sqrt{c}}-sn_{c}(2r)\right) \nonumber\\
 	&\leq \delta H_{c}(r) - \frac{2\delta ksn_{c}(r)sn_{c}'(r)}{sn_{c}^{2}(r)} + \frac{\delta k}{sn_{c}^{2}(r)}\left(\frac{2}{\sqrt{c}}-\frac{1}{\sqrt{c}}\sin(2\sqrt{c}r)\right) \nonumber \\
 	&=\delta H_{c}(r) + \frac{2\delta k}{n-1} H_{c}(r) + \frac{\delta k}{sn_{c}^{2}(r)}\left(\frac{2}{\sqrt{c}}-\frac{2\cos(\sqrt{c}r)\sin(\sqrt{c}r)}{\sqrt{c}}\right) \nonumber\\
  &=\delta H_{c}(r) + \frac{2\delta k}{n-1} H_{c}(r) + 2\delta k \left(\frac{\sqrt{c}}{\sin^{2}(\sqrt{c}r)}-\frac{\sqrt{c}\cos(\sqrt{c}r)}{\sin(\sqrt{c}r)}\right) \nonumber\\
  &=\delta H_{c}(r) + \frac{2\delta k}{n-1} H_{c}(r) +  2\delta k \left(\frac{\sqrt{c}}{\sin^{2}(\sqrt{c}r)}-\frac{sn_{c}'(r)}{sn_{c}(r)}\right) \nonumber\\
  &=\delta H_{c}(r) + \frac{2\delta k}{n-1} H_{c}(r) + 2\delta k \left(\frac{\sqrt{c}}{\sin^{2}(\sqrt{c}r)}-\frac{H_{c}(r)}{n-1}\right) \nonumber\\
  &=\left( H_{c}(r) + \frac{2 k\sqrt{c}}{\sin^{2}(\sqrt{c}r)}\right)\delta.
 \end{align}

 Developing the second term of the above expression, we obtain  
  \begin{align*}
  	H_{\eta,T}(r) &\leq \left(1+\frac{4k}{n-1}\cdot\frac{1}{\sin(2\sqrt{c}r)}\right)\delta H_{c}(r).
  \end{align*}

\end{proof}

\begin{remark}
 Note also that for the case $c=0$, we can obtain  from \eqref{tcii15} the following estimate

\begin{equation*}
		H_{\eta,T}(r)\leq \delta\left(\frac{n-1}{r}-\frac{2}{r}\eta(r)+ \frac{2}{r^2}\int_{0}^{r}\eta(t)\;dt\right).
\end{equation*}

 When $T=I$ we recover the estimate for the drifted Laplacian operator given in \cite{fzl}.  

\end{remark}

\section{Applications}\label{Application}

We need some definitions to proceed with the applications of the generalized weighted mean curvature comparison theorem. A normal minimal geodesic $\gamma:[0,\infty)\to M$ is a ray if $d(\gamma(t),\gamma(s))=|t-s|$ for all $t,s\geq 0$. Similarly, a normal minimal geodesic $\gamma:(-\infty,+\infty)\to M$ is a line if $d(\gamma(t),\gamma(s))=|t-s|$ for all $t,s\in\mathbb{R}$.
Given a ray $\gamma$, we say that $\overline{\gamma}:[0,\infty)\to M$ is an asymptote from $p$ to $\gamma$ if it is a ray which arises as limit of a sequence of minimal geodesic segments, $\sigma_{i}$, from $p$ to $\gamma(t_{i}),t_{i}\to\infty$. Asymptotes are not necessarily unique. Now, let $b^{\gamma}_{t}:M\to \mathbb{R}$, $b^{\gamma}_{t}(x)= t-d(x,\gamma (t))$ for $t\geq 0$. The family of functions $\{b^{\gamma}_{t}\}_{t\geq0}$ forms a pointwise bounded equicontinuous family that is also pointwise nondecreasing on $t$.  Therefore, the Busemann function associated to $\gamma$, defined as
	\begin{equation}
	b^{\gamma}(x)=\lim_{t\to\infty} (t-d(x,\gamma (t))), \quad x\in M,
	\end{equation}
 exists, is Lipschitz continuous on $M$ and is bounded by $d(x,\gamma(0))$. See~\cite{hb} for more details about the Busemann functions. We will deal with the nondifferentiability of Busemann functions through barrier functions. 

A lower barrier (respectively, upper barrier) for a continuous function $f$ at the point $q$ is a $C^{2}$ function $h$, defined in a neighborhood $U$ of $q$ such that $h(q)=f(q)$ and $h(x)\leq f(x)$  (respectively, $h(x)\geq f(x)$) on $U$. If $f$ is continuous, we say that $\square_{\eta}f\geq a$ at $q$ in the barrier sense if, for any $\varepsilon >0$, there is a lower barrier function $h_{\varepsilon}$ for $f$ at $q$ such that $\square_{\eta}h_{\varepsilon}\geq a-\varepsilon$. An upper bound on $\square_{\eta}$ is defined analogously in terms of upper barriers.

\begin{lem}\label{bsub} If $M$ is a complete, noncompact Riemannian manifold with $R_{\eta,T}\geq 0$ for some bounded function $\eta$ then $\square_{\eta}b^{\gamma}\geq 0$ in the barrier sense.	
\begin{proof}
	Let $q\in M$ be fixed, $b^{\gamma}$ the Busemann function associated to  $\gamma$ and let $\overline{\gamma}$ the asymptote of $\gamma$ such that $\overline{\gamma}(0)=q$ and $\overline{\gamma}'(0)=v_{0}$. For any $t>0$, we define the function
	$$h_{t}(x)=t-d(x,\overline{\gamma}(t))+b^{\gamma}(q).$$
	
	By the same arguments of~\cite{ww}, $h_{t}$ is a lower barrier for $b^{\gamma}$. Then, by \eqref{tc2ii},
\begin{equation*}
	\square_{\eta}(h_{t})(x)=\square_{\eta}(-d(x,\overline{\gamma}(t)))\geq \delta \frac{(n+4k+1)}{t}.
\end{equation*}

Taking $t\to\infty$ we obtain the desired result.
\end{proof}
\end{lem}

Since $\square_{\eta}$ is a perturbation of the Cheng-Yau operator, by a similar analysis as in~\cite{Ch,pp} it is possible to prove the following two results.

\begin{theorem}[Maximum Principle]\label{pmf} Let $M$ be a connected Riemannian manifold and $f\in C^{0}(M)$. Suppose that $\square_{\eta} f\geq 0$ in the barrier sense. Then $f$ attains no weak local maximum unless it is a constant function. In particular, if $f$ has a global maximum, then $f$ is constant.
\end{theorem}
\begin{theorem}[Regularity]\label{rge} If $f:(M,\langle, \rangle)\to \mathbb{R}$ is continuous and $\square_{\eta}f=0$ in the barrier sense, then $f$ is smooth.
\end{theorem}

We are finally ready to prove the applications.

\begin{proof}[\textbf{Proof of Theorem~\ref{split}}]
 Denote by $\gamma_{+}$ and $\gamma_{-}$ the two rays which form the line $\gamma$ and let $b^{+}$, $b^{-}$ be their Busemann functions. It is easy to verify that $b^{+} + b^{-}\leq 0$ on $M$ and $b^{+} + b^{-}=0$ on $\gamma$, and so $b^{+} + b^{-}$ has a maximum on $\gamma(0)$. By Lemma~\ref{bsub}, $\square_{\eta}b^{\pm}\geq 0$ in the barrier sense, thus $\square_{\eta}(b^{+} + b^{-})\geq 0$. From Theorem~\ref{pmf}, we see that $b^{+} + b^{-}$ is constant and equal to zero, whence $b^{+} =- b^{-}$. Therefore, $\square_{\eta}b^{+}=\square_{\eta}b^{-}=0$, which implies that $b^{+}$ and $b^{-}$ are both smooth functions by Theorem~\ref{rge}. A simple calculation using the triangle inequality shows that the asymptotes $\overline{\gamma}_{+}$ and $\overline{\gamma}_{-}$ form a line in $M$ and, when this occurs, $\overline{b}^{+}$ and ${b}^{+}$ differ from each other only by a constant. From the Bochner type formula~\eqref{forboch}, applied to the function $b^{+}$ with $|\nabla b^{+}|=1$, we obtain
\begin{equation*}\label{fboch}
	0=\nabla b^{+}(\square_{\eta}b^{+}) + R_{\eta,T} (\nabla b^{+}, \nabla b^{+}) + \langle\nabla^{2}b^{+}, \nabla^{2}b^{+}\circ T\rangle  - \langle\nabla^{2}b^{+},\nabla_{\nabla b^{+}}T\rangle.
\end{equation*}

Since $T$ is a radially parallel tensor, $\square_{\eta}b^{+}=0$ and $R_{\eta,T}\geq 0$, we have

\begin{equation*}\label{fboch2}
	0\geq tr((\nabla^{2}b^{+})^{2}\circ T) .
\end{equation*}

Therefore, $tr((\nabla^{2}b^{+})^{2}\circ T) =0$ and $(\nabla^{2}b^{+})^{2}=0$, from which it follows that $\nabla b^{+}$ is a parallel field. The rest of the proof is the same as in~\cite{ww}. See also~\cite{fzl,zhu}.
\end{proof}

\begin{proof}[\textbf{Proof of Theorem~\ref{Myers}}]
Let $p_{1}$ and $p_{2}$ be two points in $M$ with $d(p_{1},p_{2})\geq\frac{\pi}{\sqrt{c}}$ and set
\begin{equation*}
    B\coloneqq d(p_{1},p_{2}) - \frac{\pi}{\sqrt{c}}.
\end{equation*}

Let $r_{1}(x)=d(p_{1},x)$, $r_{2}(x)=d(p_{2},x)$ and denote by $e$ the excess function for the points $p_{1}$ and $p_{2}$, i.e., 
\begin{equation*}
    e(x)\coloneqq d(p_{1},x) + d(p_{2},x) - d(p_{1},p_{2}),
\end{equation*}
which measures how much the triangle inequality fails to be an equality. By the triangle inequality, $e(x)\geq 0$ and $e(\gamma(t))=0$, where $\gamma$ is a minimal geodesic segment from $p_{1}$ to $p_{2}$. Therefore, $\square_{\eta}(e)(\gamma(t))\geq 0$.
Let $y_{1}=\gamma\left(\frac{\pi}{2\sqrt{c}}\right)$ and $y_{2}=\gamma\left(\frac{\pi}{2\sqrt{c}}+ B\right)$, then $r_{i}(y_{i})=\frac{\pi}{2\sqrt{c}}$, $i=1,2$, and by~\eqref{est1} we have
\begin{equation}\label{est2}
 \square_{\eta}(r_{i}(y_{i}))\leq 2\delta k\sqrt{c}.   
\end{equation}

Since we cannot give an estimate for $\square_{\eta}(r_{1}(y_{2}))$ by directly using Theorem~\ref{tcii}, by applying Theorem~\ref{tci} and Equation~\eqref{est2} we have
\begin{equation*}
 \square_{\eta}(r_{1}(y_{2}))\leq 2\delta k\sqrt{c} - B(n-1)\delta c. 
\end{equation*}

Therefore,
\begin{equation*}
  0\leq\square_{\eta}(e)(y_{2}) = \square_{\eta}(r_{1})(y_{2}) + \square_{\eta}(r_{2})(y_{2}) 
  \leq 4\delta k\sqrt{c}- B(n-1)\delta c,
\end{equation*}
which implies $B\leq \frac{4k}{(n-1)\sqrt{c}}$  and $d(p_{1},p_{2})\leq \frac{\pi}{\sqrt{c}}+\frac{4k}{(n-1)\sqrt{c}}$.
\end{proof}

For completeness we give a proof of the maximum principle.

\begin{proof}[\textbf{Proof of Theorem~\ref{pmf}(Maximum principle)}]

It is sufficient to prove the theorem for the Cheng-Yau operator. Assume that $f$ reaches a maximum at $p\in M$, i.e., $f(p)\geq f(q)$, for all $q$ near $p$. Consider $f|_{B(p,\delta)}:B(p,\delta)\to \mathbb{R}$ the restriction of $f$ to a geodesic ball of radius  $\delta >0$ centered at $p$ with $\delta<inj(p)$, $\square f\geq 0$ in the barrier sense, and a global maximum at $p$.

If $f$ is constant there is nothing to prove, otherwise, we can assume that there exists a $q_{0}\in\partial B(p,\delta)$ such that $f(p)>f(q_{0})$. By the continuity of $f$, $f(p)>f(q)$, for all $q$ in a neighborhood of $q_{0}$ on $\partial B(p,\delta)$. Choose a normal coordinate system $(x_{1},\cdots,x_{n})$ centered at $p$ and define the smooth function
	
	\begin{equation*}
		\psi(x)=x_{1}-k(x^{2}_{2}+\cdots x^{2}_{n}),
	\end{equation*}
where $k\in\mathbb{R}$ is so large that if $y\in W=\{x\in\partial B(p,\delta); f(x)=f(p)\}$, then $\psi(y)<0$.

Note that 
\begin{itemize}
	\item [$(i)$]$W\subsetneq \partial B(p,\delta)$;
	\item [$(ii)$] $\nabla \psi\neq 0$ on $\overline{B(p,\delta)}$, since $\frac{\partial\psi}{\partial x_{1}}=1$;
	\item [$(iii)$]  $\psi(p)=0$;
	\item [$(iv)$] $ \psi<0$ on $W$.
\end{itemize}

In this way, we can construct a smooth function $h\coloneqq e^{a\psi}-1$ such that
 $h(p)=0$ and
	 $\square h>0$ on $\overline{B(p,\delta)}$. In fact, since 	
	\begin{equation*}
		\square h = div(T(\nabla h)) =div (T(\nabla (e^{a\psi}-1)))= div(T(ae^{a\psi}\nabla \psi))
	\end{equation*}
and from
\begin{equation*}
	div(T(\varphi Z))=\varphi\langle divT, Z\rangle + \varphi\langle\nabla Z,T\rangle + T(\nabla \varphi,Z),
\end{equation*}
with $\varphi \in C^{\infty}(M)$ and $Z\in\mathfrak{X}(M)$, we have
\begin{align*}
\square h &= ae^{a\psi}\langle\nabla^{2}\psi,T\rangle  + T(\nabla(ae^{a\psi}),\nabla\psi) \\
&= ae^{a\psi}\square\psi + T(a^{2}e^{a\psi}\nabla\psi,\nabla\psi) \\
&= ae^{a\psi}\square\psi + a^{2}e^{a\psi} T(\nabla\psi,\nabla\psi) \\
&\geq ae^{a\psi}\square\psi + a^{2}e^{a\psi}\varepsilon|\nabla\psi|^{2},
\end{align*}
where the last inequality comes from the fact that $\varepsilon|X|^{2}\leq \langle T(X),X\rangle$. The result follows taking $a$ sufficiently large.

For $\eta$ sufficiently small, consider the function
\begin{equation*}
	\overline{f}=f+\eta h \quad \mbox{on}\quad \overline{B(p,\delta)}.
\end{equation*}

Since $f(p)>f(x),\forall x\in \partial B(p,\delta)$, we have
\begin{equation*}
	f(p)>f(x)+\eta h(x), \forall x \in \partial B(p,\delta).
\end{equation*}

Thus, $f(p)>max\{\overline{f}(x);x\in\partial B(p,\delta)\}$. Since $\overline{f}(p) = f(p)$,  it follows that

\begin{equation*}
\overline{f}|_{\partial B(p,\delta)}< \overline{f}(p),
\end{equation*}
that is, $\overline{f}$ has a maximum at some point $z$ inside $B(p,\delta)$ provided that $\eta$ is sufficiently small.

Given that $\square f\geq 0$ in the barrier sense, if $f_{z,\varepsilon}$ is a barrier function for $f$ at $z$, then $\overline{f}_{z,\varepsilon}=f_{z,\varepsilon}+\eta h$ is also a barrier function for $\overline{f}= f+\eta h$ at $z$ and $\square f_{z,\varepsilon}\geq -\varepsilon$.
For $\varepsilon>0$ sufficiently small, we have
\begin{equation}\label{pm1}
\square \overline{f}_{z,\varepsilon}= \square (f_{z,\varepsilon}+\eta h)\geq -\varepsilon+\eta\square h>0.
\end{equation}

Since
\begin{equation*}
\overline{f}_{z,\varepsilon}(z)=\overline{f} (z)\quad\mbox{and} \quad \overline{f}_{z,\varepsilon}(x)\leq \overline{f}(x) \; (x \; \mbox{in a neighborhood of}\; z),
\end{equation*}
we have that $\overline{f}_{z,\varepsilon}=f_{z,\varepsilon}+\eta h$ also has a maximum at $z$, which implies that $\nabla^{2}\overline{f}_{z,\varepsilon}\leq 0$. Thus, $\square\overline{f}_{z,\varepsilon}=\langle\nabla^{2}\overline{f}_{z,\varepsilon},T\rangle\leq 0$, which contradicts~\eqref{pm1}. Therefore, for every sufficiently small  $\delta$, $f|_{\partial B(p,\delta)}=f(p)$. Since $M$ is connected this implies $f\equiv f(p)$.

\end{proof}

\begin{remark}
For regularity, it is sufficient to note that it is a local property, which means that this theorem is obtained directly from the regularity theory for elliptic PDEs. 
\end{remark}

\bmhead{Acknowledgements} We would like to thank the Professor Dragomir Mitkov Tsonev for useful comments.

\section*{Declarations}
\begin{itemize}
\item {\bf Funding:} The authors have been partially supported by Fundação de Amparo à Pesquisa do Estado do Amazonas (FAPEAM) and Coordenação de Aperfeiçoamento de Pessoal de Nível Superior (CAPES), Grant 001.

\item {\bf Conflict of interest/Competing interests:} the authors state that there is no conflict of interest.

\item {\bf Ethics approval and consent to participate:} not applicable.

\item {\bf Consent for publication:} open acess.

\item {\bf Data availability:} not applicable.

\item {\bf Materials availability:} not applicable.

\item {\bf Code availability:} not applicable.

\item {\bf Author contribution:} not applicable.

\end{itemize}

\end{sloppypar}

\end{document}